\documentclass[10pt]{article}

\usepackage[english]{babel}
\usepackage{amsmath,amssymb,epsfig,bbm}
\usepackage{psfrag}

\usepackage{color}

\DeclareMathOperator{\supp}{supp}

\begin{document}

\numberwithin{equation}{section}

\newtheorem{thm}{Theorem}
\newtheorem{defi}{Definition}
\newtheorem{lem}{Lemma}

\newcommand{\pf}{\noindent{\bf Proof: }}
\newcommand{\qed}{ \hspace*{3em} \hfill $\square$}

\newcommand{\N}{\mathbbm{N}}
\newcommand{\Z}{\mathbbm{Z}}
\newcommand{\R}{\mathbbm{R}}
\newcommand{\pr}{\mathbbm{P}}
\newcommand{\ex}{\mathbbm{E}}

\newcommand{\Bo}{\mathcal{B}}
\newcommand{\F}{\mathcal{F}}
\newcommand{\Dset}{\mathcal{D}}
\newcommand{\X}{\mathcal{X}}
\newcommand{\Aset}{\mathcal{A}}

\newcommand{\de}{\delta}
\newcommand{\ep}{\epsilon}
\newcommand{\be}{\beta}
\newcommand{\ga}{\gamma}
\newcommand{\Ga}{\Gamma}
\newcommand{\Gac}{\Ga_{\cup}}
\newcommand{\ph}{\varphi}
\newcommand{\De}{\Delta}
\newcommand{\la}{\lambda}
\newcommand{\si}{\sigma}

\newcommand{\gro}{G}
\newcommand{\esc}{L}
\newcommand{\rp}{d_0}
\newcommand{\die}{d_1}
\newcommand{\dio}{d_2}

\newcommand{\sli}{l}

\newcommand{\leo}{\zeta}
\newcommand{\tleo}{\tilde{\zeta}}

\newcommand{\Mli}{M_{\text{in,fi}}}
\newcommand{\Mpl}{M'_{\text{in}}}
\newcommand{\Mi}{M_{\text{fi}}}
\newcommand{\Mpi}{M'_{\text{fi}}}
\newcommand{\Mpli}{M'_{\text{in,fi}}}
\newcommand{\Mpdel}{M'_{\text{de}}}

\newcommand{\Din}{\De_{\text{in}}}
\newcommand{\Dfi}{\De_{\text{fi}}}
\newcommand{\Dde}{\De_{\text{de}}}

\newcommand{\Zc}{Z^{\cup}}
\newcommand{\tZ}{\tilde{Z}}
\newcommand{\tZc}{\tZ^{\cup}}
\newcommand{\tze}{\tilde{\zeta}}
\newcommand{\tzed}{\tze'}
\newcommand{\tZd}{\tZ'}
\newcommand{\tiC}{\tilde{C}}
\newcommand{\tGa}{\tilde{\Gamma}}
\newcommand{\tGac}{\tGa_{\cup}}
\newcommand{\ttau}{\tilde{\tau}}

\newcommand{\CRM}{CR-model}
\newcommand{\CRMs}{CR-models}

\newcommand{\eff}{\text{eff}}
\newcommand{\peff}{\text{peff}}
\newcommand{\normal}{standard }
\newcommand{\Normal}{Standard }
\newcommand{\bem}[1]{\textcolor{blue}{\emph{#1}}} 

\begin{center}
{\bf \Large The 2-type continuum Richardson model:\\ 
Non-dependence of the survival of both types\\
 on the initial configuration}

\bigskip 

Sebastian Carstens and Thomas Richthammer 
\footnote[1]{Mathematisches Institut der LMU, 
Theresienstr. 39, 80333 Muenchen, Germany\\
Research supported in part by NSF grant DMS-0300672}
\end{center}


\begin{abstract}
We consider the model of Deijfen et al.\ for competing growth of 
two infection types in $\R^d$, based on the Richardson model on $\Z^d$.
Stochastic ball-shaped infection outbursts transmit the infection type 
of the center to all points of the ball 
that are not yet infected. Relevant parameters of the model 
are the initial infection configuration, the (type-dependent) growth 
rates and the radius distribution of the infection outbursts. 
The main question is that of coexistence: 
Which values of the parameters allow the unbounded growth 
of both types with positive probability? 
Deijfen et al.\ conjectured that the initial configuration 
basically is irrelevant for this question, and gave a proof for this 
under strong assumptions on the radius distribution, 
which e.g.\ do not include the case of a deterministic radius. 
Here we give a proof that doesn't rely on these assumptions. 
One of the tools to be used is 
a slight generalization of the model with 
immune regions and delayed initial infection configurations.\\

\noindent
Key words: Continuum growth model, Richardson's model, initial configuration, 
competing growth, shape theorem.\\

\noindent 
Classification: Primary 60K35, Secondary 82B43.
\end{abstract}



\section{Introduction} \label{secintro}

We consider a model for the competing growth of $n$ types of infections on $\R^d$, 
where  $d \ge 1$ is the number of spatial dimensions. 
We will refer to it as the ($d$-dimensional $n$-type) continuum Richardson model (\CRM).
The corresponding model with one type was introduced by M.~Deijfen in \cite{D}
as a continuum version of the growth model on $\Z^d$ introduced by 
D.~Richardson in \cite{R}. The multitype versions of these models 
were first considered by O.~H\"aggstr\"om and 
R.~Pemantle in \cite{HP} (discrete model) and 
M.~Deijfen, O.~H\"aggstr\"om and J.~Bagley in \cite{DHB} (continuum model).\\

The \CRM\ is a stochastic process $Z = (Z^i_t)_{i,t}$, where $Z^i_t$ denotes 
the subset of $\R^d$ that is infected with type $i \in \{1,\ldots,n\}$ 
(i.e.\ $i$-infected) at time $t \ge 0$. Initially, 
given disjoint regions $\Ga_i$ are $i$-infected.  
Whenever a region is $i$-infected, it stays like that.  
Furthermore it tries to $i$-infect healthy "neighboring" regions 
by means of stochastic infection outbursts. 
The waiting time for the first outburst of infection type $i$ after time $t$ 
is exponentially distributed with rate $\be_i\la^d(Z_t^i)$, 
where $\la^d$ denotes Lebesgue measure and $\be_i$ is a type-dependent growth rate. 
The outburst has the shape of a ball, where the center is chosen uniformly in $Z_t^i$ 
and the radius is chosen w.r.t.\ a given radius distribution $\rho$. At the time 
of the outburst all points within the ball that are not yet infected get $i$-infected.
We note that the $i$-infected region $Z^i_t$ is increasing in $t$, 
and at any given time $t$ the sets $(Z^i_t)_i$ are disjoint.\\

The parameters of the \CRM{} are the initial configuration $\Ga = (\Ga_1,\ldots,\Ga_n)$, 
the growth rates $\be = (\be_1,\ldots,\be_n)$ and the radius distribution $\rho$. 
Throughout the paper we will assume that $\la^d(\Ga_i) > 0$, $\be_i > 0$ 
and $\rho(\{0\}) = 0$ to avoid trivialities; furthermore we will assume that 
$\Ga_i$ is bounded and $\rho$ satisfies
\begin{equation} 
\int_0^{\infty} r^d d\rho(r) \, < \, \infty \label{momentr}, 
\end{equation}
which will ensure that the model is well defined. Other assumptions on $\rho$ 
that sometimes play an important role are stronger moment assumptions such as 
\begin{equation} 
\int_0^\infty e^{\ph r} d\rho(r)  \, < \, \infty \quad 
\text{ for some } \quad  \ph > 0,\label{mogenr}
\end{equation}
and conditions on the support of $\rho$ such as 
\begin{equation} \label{esssup0}
\rho((0,\de)) >0 \quad \text{ for all } \quad \de >0. 
\end{equation}
For instance, \cite{D} and \cite{DHB} show for the 1-type \CRM\ 
that under condition  \eqref{mogenr} on $\rho$ 
the asymptotic shape of $Z_t/t$ for $t \to \infty$ is a ball; 
we have stated this result as Theorem~\ref{thmshape} 
in Section~\ref{secpropmod}. \\

An interesting question for the \CRM\ is that of coexistence: 
Do all infection types grow unboundedly at the same time with positive probability? 
In the 2-type \CRM\ it is conjectured
that we have coexistence iff $\be_1 = \be_2$. It is known that 
for fixed $\be_1$ at most countably many values of $\be_2$ allow coexistence
(see M.~Deijfen, O.~H\"aggstr\"om and J.~Bagley in \cite{DHB}) 
and the value $\be_2 = \be_1$ is one of them 
(see M.~Deijfen and O.~H\"aggstr\"om in \cite{DH1}).
Both results assume condition \eqref{mogenr} 
and concern models with special initial configurations 
consisting of two disjoint balls. However, the question of 
coexistence is basically independent of the initial configuration. 
This is shown in \cite{DHB} assuming that $\rho$ satisfies 
\eqref{mogenr} and \eqref{esssup0}. 
The aim of the present article is to show how these extra assumptions can be avoided. 
(In fact it turns out that assuming condition \eqref{momentr} 
is sufficient.) 
As an immediate consequence the aforementioned coexistence results 
extend to basically all initial configurations for all 
radius distributions satisfying \eqref{mogenr}.
This is desirable as some of the most natural 
choices of $\rho$ do not satisfy  \eqref{esssup0}, e.g.\ 
the case of outbursts with a deterministic radius.
\\

For stating our result let us consider a 2-type \CRM\ 
with initial condition $\Ga = (\Ga_1,\Ga_2)$ and growth rates $\be = (\be_1,\be_2)$. 
Let $B_{\Ga}$ denote the smallest ball centered at the origin containing 
$\Gac := \Ga_1 \cup \Ga_2$. Let 
$\esc_i$ be the event that type $i$ leaves $B_{\Ga}$,
$\gro_i$ the event that type $i$ reaches points 
arbitrarily far from the origin, 
and  $\gro := \gro_1 \cap \gro_2$ 
the event of unbounded growth of both types. Whenever we consider more 
than one model, we will indicate the relevant parameters in 
parentheses after the corresponding event.
\begin{thm} \label{thmgro}
We consider two $d$-dimensional 2-type \CRMs\ 
with $d \ge 2$, initial configuration $\Ga$ and $\Ga'$ respectively, 
growth rates $\be_i = \be_i'$ and radius distributions 
$\rho = \rho'$ (satisfying \eqref{momentr}). 
If we have $\pr(\esc_i(\Ga)) > 0$ and 
$\pr(\esc_i(\Ga')) > 0$ for $i=1,2$, 
the possibility of coexistence doesn't depend on the 
initial condition: 
\[
\pr(\gro(\Ga)) >0  \quad \iff \quad \pr(\gro(\Ga')) > 0.
\]
\end{thm}
We note that the condition on the events  
$\esc_i$ merely states that in the initial configurations no type strangles the other. 
As mentioned above, combining Theorem~\ref{thmgro} with the coexistence results 
from \cite{DH1} and \cite{DHB} immediately gives the following: 
\begin{thm} \label{thmcoex}
We consider a $d$-dimensional  2-type \CRM\ 
with $d \ge 2$, initial configuration $\Ga$, growth rates $\be = (\be_1,\be_2)$
and radius distribution $\rho$ satisfying \eqref{mogenr}.
\begin{enumerate}
\item[(a)] If $\be_1 = \be_2>0$ and  $\pr(\esc_i(\be)) > 0$ for $i=1,2$
we have $\pr(\gro(\be)) >0$.
\item[(b)] For fixed $\be_1$ we have $\pr(\gro(\be))=0$ 
for all but countably many values of $\be_2$.
\end{enumerate}
\end{thm}

We briefly discuss to what extent the proof of Theorem~\ref{thmgro} in \cite{DHB} 
makes use of the stronger assumptions on the radius distribution. 
One of the tools used is the aforementioned theorem on the asymptotic shape 
of the infected region which relies on \eqref{mogenr}. 
But only the lower bound of the asymptotic shape is used, 
which is still valid if the radii of the outbursts are increased. 
Thus it is easy to relax \eqref{mogenr} to \eqref{momentr}. 
Another tool in the proof is the construction of an infection evolution 
that infects some given points with type 1 (``type-1 points'')
and some other given points with type 2 (``type-2 points''). 
For this, condition \eqref{esssup0} is essential: 
Suppose that some type-1 points are surrounded by type-2 points. 
Utilizing \eqref{esssup0} it is easy to construct a sequence of sufficiently small outbursts 
that 1-infect the type-1 points without 1-infecting the type-2 points. 
Without \eqref{esssup0} this seems to be hopeless. 
Instead, our strategy will be to 2-infect all type-2 points, 
but to 1-infect only a single (suitably chosen) type-1 point. 
The construction of such an infection evolution turns out to be possible without \eqref{esssup0}, 
though somewhat complicated in terms of geometry. 
We then have to investigate how the unknown infection states of the remaining type-1 points 
affect the subsequent infection evolution. 
Our key tool for this is a generalized version of the \CRM\ 
that allows immune regions and delayed infections in its initial configuration. \\

We define a generalized initial configuration $\Ga$ of an $n$-type \CRM\ 
to be a finite collection $\Ga = (\Ga_j)$ of disjoint bounded Borel subsets 
of $\R^d$, each of which has an associated type $i(\Ga_j) \in \{1,\ldots,n\}$ 
and an associated time $t(\Ga_j) \in [0,\infty]$. $\Ga_j$ is considered to be 
uninfected at times $t < t(\Ga_j)$ and infected with type $i(\Ga_j)$ at times 
$t \ge t(\Ga_j)$. If $t(\Ga_j) < \infty$ this corresponds to a delayed 
initial configuration, and if $t(\Ga_j) = \infty$ this corresponds to an immune 
region (and in this case $i(\Ga_j)$ is irrelevant). 
We will use the shorthand notation $\Gac := \bigcup_{j} \Ga_{j}$. 
Additionally we are given growth parameters $\be_i > 0$ 
and a radius distribution $\rho$ as before. 
The corresponding stochastic process will again be 
denoted by $Z = (Z_t^i)_{i,t}$, 
where $Z_t^i$ denotes the $i$-infected subset of $\R^d$ at time $t \ge 0$. 
The dynamics of outbursts is a generalization of that for 
the \CRM\ with \normal initial configurations. 
It is defined by the following properties: 
\begin{itemize}
\item $Z$ is a Markov process (which is time-homogeneous iff 
$t(\Ga_j) \in \{0,\infty\}$  $\forall j$).
\item  In $\Gac$ we have a deterministic infection evolution 
as described above.
\item Regions of $\Gac^c$ can only be infected by ball shaped 
infection outbursts. If the center of such a ball is $i$-infected, 
then the outburst $i$-infects all points 
within the ball that are not yet infected and not in $\Gac$.
\item At a given time $t$ an outbursts of type $i$ occurs 
at rate $\be_i \la^d(Z^i_t)$.
The center of the corresponding ball is chosen uniformly 
in $Z_t^i$ and the radius is chosen w.r.t. $\rho$.
\end{itemize}
The above properties uniquely define the distribution of the process $Z$. 
We will refer to such a process as \CRM\ 
with generalized initial configuration or generalized \CRM.\\ 

This article is organized as follows: In Section~\ref{secset} 
we clarify some notation, describe two ways 
to construct a generalized \CRM\ from a Poisson point process, 
and state a shape theorem and some other important properties of the generalized \CRM. 
The proof of the 
corresponding lemmas and theorems is relegated to Section~\ref{secproofset}. 
In Section~\ref{secproofthm} we prove Theorem~\ref{thmgro}. 
The proof of the corresponding lemmas is relegated to Section~\ref{secprooflemthm}.


\section{Setting} \label{secset}

\subsection{The underlying space and point process}

First we introduce some notation concerning subsets of $\R^d$ and point processes. 
We denote the Lebesgue measure on $\R^d$ by $\la^d$ and the restriction of 
$\la = \la^1$ to $\R_+ = [0,\infty)$ by $\la_+$.
Let $B(x,r)$ denote the closed ball with center $x\in \R^d$ and radius $r>0$. 
Let $d(x,y) := |x-y|$ be the Euclidean distance between 
two points $x,y \in \R^d$. We will use $d(.,.)$ in the usual 
way also for the distance of a point and a set or the 
distance of two sets.
For $A \subset \R^d$ and  $r >0$ let 
\[
A_{+r} = \{x \in \R^d: d(x,A) \le r\}. 
\]

Considering point processes we will be mainly concerned with 
points $p = (x,s,r,w) \in \R^d  \times \R_+ \times \R_+ \times \R_+
=: \R^d_{\times}$. We will refer to the components of such 
a point $p$ as position $x$, time $s$, radius $r$ and strength $w$.
For $A \subset \R^d$ and time interval $I \subset \R_+$ 
by abuse of notation we will sometimes consider $A$ and $A \times I$
as subsets of $\R^d_{\times}$. 
A point configuration 
$X$ is defined as a locally finite subset of $\R^d_{\times}$.  
Let $\X$ denote the set of all point configurations. 
For a Borel set $E\subset\R^d_{\times}$ and $X \in \X$
let $N_E(X)$ be the number of points of $X$ in $E$, and let $\F_E$ be 
the $\si$-algebra on $\X$ generated by all counting variables
$N_{E'}$, where $E'$ is a Borel subset of $E$. 
$\F := \F_{\R^d_{\times}}$ will be the $\si$-algebra usually 
associated to $\X$. Sets $M \in \F_E$ will be called 
events depending only on (the information in) $E$.
For a given radius distribution $\rho$ we consider the Poisson point process with intensity measure 
$\la^d_\times := \la^d \otimes \la_+ \otimes \rho \otimes \la_+$ and 
denote its distribution by $\pr$. 


\subsection{Constructing the model from a Poisson point process}

Many proofs rely on the comparison of different generalized \CRMs, 
so we have to construct suitable couplings between them. 
To prepare this, we describe two methods to construct an infection evolution 
$Z = Z(X)$ for a given point configuration $X \in \X$. 
Typically, $X$ will be a realization of $\pr$. \\

The first method to construct an infection evolution from $X$ 
is called "scanning from time of infection". The infection 
evolution in the generalized initial configuration $\Gac$ is deterministic: 
Every region $\Ga_j$ is infected at the 
associated time with the associated type (and is uninfected before). 
Whenever a previously uninfected region $A$ is $i$-infected at some time $t$, 
it gets a scanning device that starts scanning for points of $X$ in 
$A \times [t,\infty) \times \R_+ \times [0,\be_i]$. A point $p = (x,s,r,w) \in X$
located in $A$ with time $s \ge t$ and strength $w \le \be_i$ will be scanned at time $s$. 
At this time the point $p$ produces an outburst $i$-infecting all points of $B(x,r)$ 
that are not yet infected and not in $\Gac$. 
We stress that in this construction a point $p = (x,s,r,w) \in X$ can produce an outburst 
only if $x$ has been infected by time $s$; if $x$ has been $i$-infected before time $s$ 
and $w \le \be_i$, $p$ produces an outburst at time $s$.\\

This construction is changed slightly in the second method: 
"scanning from time 0". The infection evolution in $\Gac$ 
is the same as above, and again a previously uninfected region $A$ 
that is $i$-infected at some time $t$ gets a scanning device. But now this device 
starts scanning for points of $X$ in $A \times [0,\infty) \times \R_+ \times [0,\be_i]$.
So a point $p = (x,s,r,w) \in X$ located in $A$ with time $s \ge 0$ and strength $w \le \be_i$ 
will be scanned at time $t+s$. The outburst generated by $p$ at this time has the same 
effect as above.  We stress that in this construction the only 
condition for a point $p = (x,s,r,w) \in X$ to produce an outburst is that $w$ 
is sufficiently small (depending on the eventual infection type of $x$); 
if $x$ gets $i$-infected at time $t$ and $w \le \be_i$, $p$ 
produces an outburst at time $t+s$.\\

If $X$ is chosen according to $\pr$, 
both constructions give a random infection evolution with the properties 
described in Section~\ref{secintro}, and thus yield a generalized \CRM.
Furthermore for both methods of construction there are only finitely many outbursts  
in any finite amount of time a.s.. This can be shown as in \cite{D}, Prop.~2.1, 
for the usual 1-type \CRM, making use of property \eqref{momentr} of the 
radius distribution $\rho$. As a consequence of the properties of the Poisson point process, 
in both methods of constructions a.s.\ no point is scanned at the same time
as another one or at a time associated to part of the initial configuration.
The only instance where we will use "scanning from time 0" 
is in the proof of Theorem~\ref{thmshapegen}. In every other 
instance when a \CRM\ is related to a point process or 
a point configuration we will assume that the \CRM\ is constructed 
by "scanning from time of infection".
Whenever we consider two or more initial configurations we write $Z^i_t(\Ga)$ to indicate which initial configuration the process uses. 
In the special case of a \normal initial configuration, 
the above constructions are slightly easier, see \cite{D}, 
\cite{DHB} and \cite{DH1}.\\

We conclude this subsection by introducing some notions 
useful for describing an infection evolution. 
The total infected region at time $t$ will be denoted by 
$\Zc_t := Z_t^1 \cup \ldots \cup Z_t^n$. 
Next we consider a point $p = (x,s,r,w)$  scanned at time $t$. 
$p$ is said to produce an effective outburst 
if it produces an outburst that increases the infected region, i.e. 
$B(x,r) \not \subset \Zc_t \cup \Gac$. 
$p$ is said to produce a virtual effective outburst 
if $w \le \max_i \be_i$ and $B(x,r) \not \subset \Zc_t \cup \Gac$ 
(but $x$ is not necessarily infected at time $s$, 
so $p$ might not produce any outburst). 
A virtual effective outburst is said 
to virtually infect points of a set $C$ if $B(x,r) \cap C \not \subset \Zc_t \cup \Gac$.
A time $t$ is called a time of growth if at this time the infected region increases, 
which may be due to the initial configuration or an effective outburst. 
Many proofs of assertions about an infection evolution are by induction on times of growth
(considered in linear order).
Finally we call a sequence $(p_n)_{n}$ of points 
$p_n$ an $i$-infection path if 
every $p_n$ produces an effective outburst that $i$-infects
the position of $p_{n+1}$.


\subsection{Auxiliary results} \label{secgrowth} \label{secpropmod}

Let $Z$ denote a 2-type \CRM\ with 
\normal initial configuration $\Ga = (\Ga_1,\Ga_2)$. In the proof of 
Theorem \ref{thmgro} we have to construct various infection evolutions. 
A building block of these constructions is to let one type, say type 1, 
grow along a given set $C$: If at time $T_1$ 
a neighborhood of some point $c_0 \in C$ is 1-infected 
and $C$ is not 2-infected, we would like to construct 
an infection evolution in the time interval $(T_1,T_2]$ 
so that type 2 doesn't grow at all and type 1 infects some enlargement 
of $C$ (as far as that is possible). For this we would like to use outbursts 
of a certain size, so we choose a $\rp >0$ with 
\begin{equation} \label{radprorp}
\rho([\rp,\rp+\ep])>0 \quad \text{ for all } \ep >0, 
\end{equation}
e.g.\ $\rp := \inf\{r > r_0: \rho([0,r]) > 
\rho([0,r_0])\}$ for some small $r_0$ with $\rho([0,r_0])<1$. 
We fix the value of such a $\rp$ for the rest of the article. 
Definition~\ref{defigrowth} gives a more precise description of the growth of type 1 
along a set $C$, and in Lemma~\ref{legrowth} we give conditions on $C$ under which 
this kind of growth is possible.
\begin{defi} \label{defigrowth} \emph{(Growth along a set)}
Let $0 < \de < {\rp}/{2}$, $0 \le T_1 < T_2$, $B,C \subset \R^d$ bounded Borel sets 
with $C_{+\rp+2\de} \subset B$ and $c_0 \in C$.  
An event $M$ of point configurations is said to describe the growth of type 1 along $C$ 
(within $B$, starting at $c_0$) in the time interval $(T_1,T_2]$ 
using outbursts of size $\rp$ with precision $\de$ if it has the following properties: 
$M$ only depends on $B \times (T_1,T_2]$, $\pr(M)>0$ and 
\begin{equation} \label{growthalong}
\begin{split}
&\{B(c_0,\de) \subset Z_{T_1}^1, \;  \Zc_{T_1}\subset B, \;  C_{+\de} \cap Z_{T_1}^2 
= \emptyset\} \; \cap \; M \\
&\subset  \quad \{ C_{+\rp-2\de} - Z^2_{T_1} \subset Z^1_{T_2}   
\subset C_{+\rp+2\de} \cup Z^1_{T_1} , \; 
Z^2_{T_2} = Z^2_{T_1}\}.
\end{split}
\end{equation}
\end{defi}
\begin{lem} \label{legrowth} \emph{(Constructing infection evolutions)} 
In the situation of Definition \ref{defigrowth} suppose that the set $C$ 
is of one of the following types: 
\begin{enumerate}
\item[(a)]
$C = \{c_0,\ldots,c_n\}$, where $d(\{c_0,\ldots,c_{i-1}\},c_i)\le \rp - 2 \de$ for all 
$i \ge 1$.
\item[(b)] $C$ is a bounded, pathwise connected set.
\end{enumerate}
Then it is possible to let type 1 grow along $C$ 
(within $B$, starting at $c_0$) in the time interval $(T_1,T_2]$ 
using outbursts of size $\rp$ with precision $\de$.
\end{lem}
We note that Lemma~\ref{legrowth} and Definition~\ref{defigrowth} apply analogously 
when the roles of type 1 and 2 are interchanged. 
An important special case of Lemma~\ref{legrowth}, case (b), 
is that $C$ is the trace of some continuous curve 
$\ga :[0,1] \to \R^d$ with $c_0 =\ga(0)$. \\

For given \CRM\ and a Borel set $B$,  let $\tau_B$ denote the first time 
all points of $B$ are infected, 
and $\leo_B$ denote the time of the last effective outburst 
produced by a point with position in $B$. These random variables have finite values: 
\begin{lem} \label{leeffective} 
Let $Z$ be a $d$-dimensional $n$-type \CRM\ with \normal initial 
configuration, and let $B \subset \R^d$ be a bounded Borel set. We have a.s.: 
\begin{enumerate}
\item[(a)]
$\tau_B < \infty$, 
i.e.\ $B$ is completely infected after a finite time.
\item[(b)]
$\leo_B < \infty$, i.e.\ only finitely many effective outbursts 
originate from $B$.
\end{enumerate}
\end{lem}
Very similar results have been obtained before. Note that in 
(b) we are not assuming condition \eqref{mogenr} 
on the radius distribution (compare to Lemma 4.5(a) in 
\cite{DHB}). The proof of part (b) basically relies on the
shape theorem stated below. For its proof we refer to \cite{DHB}. 
\begin{thm} \label{thmshape}
Let $Z$ be a $d$-dimensional 1-type \CRM\
with \normal initial configuration $\Ga$, growth rate $\be$ 
and radius distribution $\rho$ satisfying \eqref{mogenr}.
There is a $\mu>0$ (independent of $\be, \Ga$) 
such that for all $0<\ep<1$ we have a.s.: 
\[
(1-\ep)B(0,\be / \mu) \subset Z_t/t \subset
(1+\ep)B(0,\be / \mu) \quad 
\text{ for sufficiently large $t$}.
\]
\end{thm}

We need a version of Lemma~\ref{leeffective} 
for a \CRM\  $\tZ$ with generalized initial configuration $\tGa$. 
There may be regions that never will be infected 
since they are enclosed by a thick layer of an immune region. 
So we consider $B_{\tGa}$, 
the smallest ball enclosing $\tGac$ centered at the origin, 
and $L$, the event that the infection leaves $B_{\tGa}$.
Let $\tleo_B$ be the time of the last virtual outburst that is produced by a point with 
position in $B$ and virtually infecting points outside of $B_{\tGa}$.
\begin{lem} \label{leeffectivegen} 
Let $\tZ$ be a $d$-dimensional $n$-type CR-Model with 
generalized initial configuration $\tGa$, and let $B \subset \R^d$ be a bounded Borel set.
On $L$ we have a.s.:
\begin{enumerate}
\item[(a)]
$\tau_{B-B_{\tGa}} < \infty$, i.e.
$B-B_{\tGa}$ is completely infected after a finite time.
\item [(b)] 
$\tleo_B < \infty$, i.e. 
there are only finitely many virtual effective outbursts originating from $B$ 
virtually infecting points  outside of $B_{\tGa}$.
\end{enumerate}
\end{lem}
Again part (b) relies on a corresponding version of the shape theorem: 
\begin{thm} \label{thmshapegen}
Let $\tZ$ be a $d$-dimensional 1-type \CRM\
with generalized initial configuration $\tGa$, growth rate 
$\be$ and radius distribution $\rho$ satisfying \eqref{mogenr}.
There is a $\mu>0$ (independent of $\be, \tGa$) such that for all $0<\ep<1$  
we have a.s.\ on~$L$:
\[
(1-\ep)B(0,\be / \mu) \subset (\tZ_t \cup B_{\tGa})/t \subset
(1+\ep)B(0,\be / \mu) \quad 
\text{ for sufficiently large $t$}.
\]
\end{thm}


\section{Proof of auxiliary results} \label{secproofset}

For the following proofs we set 
$\be_{\wedge}  := \min_i \be_i$ 
and $\be_{\vee}  :=  \max_i \be_i$.

\subsection{Constructing infection evolutions: Lemma~\ref{legrowth}}

In both cases (a) and (b) we define an event $M$ that easily can be checked to 
induce the desired infection evolution.
For (a) let  $t_i := T_1 + \frac{i}{n+1}(T_2-T_1)$ ($0 \le i \le n+1$) and  
$B_i := B(c_i,\de)$ ($0 \le i \le n$). 
Let $M$ be the set of all point configurations with exactly one point 
in each of $B_i \times (t_{i},t_{i+1}] \times [\rp,\rp+\de) \times 
[0,\be_1]$ 
for all $0 \le i \le n$, and no other points apart from these in 
$B \times (T_1,T_2] \times [0,\be_{\vee}]$.

For (b) we note that the compact set 
$C_{+\rp-2\de}$ is covered by the union of all open balls with center
in $C$ and radius $\rp - \de$. We choose a finite subcover 
consisting of balls with centers $c_1,\ldots,c_n$. 
W.l.o.g.\ we may assume that $d(c_{i-1},c_i) \le \rp - 2 \de$ 
for $1 \le i \le n$ (by joining the points $c_i$  
by continuous curves $\subset C$ and adding points of these curves 
that are sufficiently close to each other 
to this list of centers). Let $M$ be chosen for 
$C'=\{c_0,\ldots,c_{n}\}$ according (a).


\subsection{\Normal initial configuration: Lemma~\ref{leeffective}}

For the proof of part (a) we may assume that 
$B$ is a ball with $B \supset \Ga$ (by replacing $B$ with a 
sufficiently big ball). Let $\de < \rp/2$. 
As in the proof of Lemma~\ref{legrowth}, case (b), 
we can choose $c_1,\ldots,c_n \in B$ such that 
$d(c_{i-1},c_i) \le \rp - 2 \de$ and $\bigcup_i B(c_i,\rp -\de) \supset B$.
We set $C_0 = \Ga$ and $C_i =B(c_i,\de)$, and for given $T>0$ we set  
$t_i := \frac{i}{n+1}T$ ($0 \le i \le n+1$). 
Let $M_T$ be the set of all point configurations with at least 
one point in each of $C_i \times (t_{i},t_{i+1}] \times [\rp,\infty)
\times [0,\be_{\wedge}]$ for all $0 \le i \le n$. On $M_T$ $~B$ is 
completely infected at time $T$, and for $T \to \infty$ we have 
$\pr(M_T) \to 1$. \\

For the proof of part (b) we first show that there is a constant $v >0$ s.t.  
\begin{equation} \label{shapecomparison}
M_T := \{ B(0,2tv) \subset \Zc_t \text{ for all  }t \ge T\} 
\uparrow \X \text{ a.s. \quad for }T \to \infty.
\end{equation}
This basically is a consequence of a comparison argument 
and the shape theorem (Theorem~\ref{thmshape}): 
Let $f:[0,\infty) \to [0,1)$ be a bijective increasing function 
such that $f(x) \le x$, e.g. $f(x) = x/(1+x)$. 
$\tilde{\rho} := \rho \circ f^{-1}$ has bounded support and 
thus satisfies \eqref{mogenr}. Let $\tZ$ be the 1-type \CRM\ with 
initial configuration $\Gac$, growth rate $\be_{\wedge}$ 
and radius distribution $\tilde{\rho}$. 
Theorem~\ref{thmshape} gives \eqref{shapecomparison} for $\tZ$ instead of $Z$. 
Thus it suffices to construct a coupling of $Z$ and $\tZ$ such that 
$\tZ_t \subset \Zc_t$ for all $t \ge 0$ a.s.. Let $Z$ be constructed from $X$, 
and $\tZ$ from $\tilde{X} := \{(x,s,f(r),w): (x,s,r,w) \in X\}$, 
scanning from time of infection, where $X$ is chosen w.r.t. $\pr$. 
Using this coupling the desired inclusion easily follows 
by induction on the times of growth of $\tZ$. This shows 
\eqref{shapecomparison}.

Suppose $T$ is chosen so large that $B \subset B(0,Tv)$. 
On $M_T$ every point $p = (x,s,r,w)$ producing an effective 
outburst in $B$ after time $T$ has the property 
$r \ge 2sv - Tv \ge sv$
and $w \le \be_{\vee}$. More formally $N^{\eff}_{B,T}  \cdot 1_{M_T} \le N_{E(B)}$
, where $N^{\eff}_{B,T}$ denotes the number of effective outbursts 
w.r.t. $Z$ in $B$ after time $T$ and 
\[
E(B) \, := \, \{(x,s,r,w) \in \R^d_{\times}: 
x \in B, r \ge sv, w \le \be_{\vee}\}.
\]
The $\pr$-expectation of $N_{E(B)}$ is 
\[ 
\la^d_{\times}(E(B)) \, 
= \, \int_B d\la^d(x) \int d\rho(r) \int_0^{\be_{\vee}} \!\!\!\! d\la^1(w) 
\int_0^{r/v} \!\!\! \!\!d\la^1(s)  \\ 
= \; \tfrac{\be_{\vee} \la^d(B)}{v} \int \rho(dr) r,  
\]
where the last term is finite by condition \eqref{momentr} on $\rho$. So 
\[
\pr(\{N^{\eff}_{B,0} = \infty \}\cap M_T) 
=  \pr(\{N^{\eff}_{B,T} = \infty \}\cap M_T) 
\le \pr(N_{E(B)} = \infty)
\; = \; 0,
\]
i.e. on $M_T$ the number of effective outbursts in $B$ is finite a.s., 
and the result follows from \eqref{shapecomparison}. 


\subsection{Generalized initial configuration: 
Lemma~\ref{leeffectivegen} }

For the proof of (a) we may assume that $B$ is some big  ball 
centered at $0$ minus the set $B_{\tGa}$. 
For rational parameters $c_0 \in B$, $t_0 >0$ and $\de <\rp/2$ let  
$L^{t_0}_{c_0,\de}$ be the event that $Z_{t_0}^{\cup} \cap B$ contains 
$B(c_0,\de)$. 
As $\bigcup_{t_0, c_0, \de} L^{t_0}_{c_0,\de}= L$ it suffices 
to show that $\tau_B < \infty$ a.s.\ on $L^{t_0}_{c_0,\de}$ 
for arbitrary $t_0$, $c_0$, $\de$ as above. 
This follows using the same construction as in the proof of 
Lemma~\ref{leeffective}, part (a), using $C_0 := B(c_0,\de)$. \\

The proof of part (b) is almost exactly the same 
as that of part (b) of Lemma~\ref{leeffective}. Here we let 
$N^{\peff}_{B,T}$ be the number of virtual effective outbursts in $B$ 
virtually infecting points outside of $B_{\tGa}$ after time $T$. 
We can use a comparison 
argument and Theorem~\ref{thmshapegen} to show that 
\[
M_T := \{ B(0,2tv) \subset \tZ_t \cup B_{\tGa} \text{ for all  }t \ge T\} 
\uparrow L \text{ a.s. \quad for }T \to \infty
\]
for some $v > 0$, and $N^{\peff}_{B,T} \cdot 1_{M_T} \le N_{E(B)}$  
gives $\pr(\{N^{\peff}_{B,0} = \infty \}\cap M_T) = 0$.


\subsection{Generalized Shape Theorem:  
Theorem~\ref{thmshapegen}}

The basic idea of the proof is to compare the \CRM\ $\tZ$ 
with generalized initial configuration $\tGa$ to the \CRM\ $Z$ 
with \normal initial configuration $\Ga := B_{\tGa}$, 
and then to use the shape theorem for $Z$. In order to 
be able to compare $\tZ$ and $Z$, we construct them 
from the same underlying Poisson point 
process ``scanning from time $0$'', see Section~\ref{secset}. 
For given rational times $0 \le T_1 \le T_2$ we define 
\[
M_{T_1,T_2} := \{ \leo_{\Ga}(Z) \le T_1, Z_{T_1} \subset \tZ_{T_2} \cup \Ga\} \cap L,
\] 
where $\leo_{\Ga}(Z)$ is the time of the last effective outburst 
in $\Ga$ w.r.t.\ $Z$. $Z_{T_1}$ is bounded and 
by Lemmas~\ref{leeffective} and  \ref{leeffectivegen} 
both $\leo_{\Ga} (Z) $ and $\tau_{B - \Ga} (\tZ )$ are finite on $L$ 
for any bounded Borel set $B$, where 
$\tau_{B - \Ga} (\tZ )$ is the first time that $B$ is infected w.r.t.\ $\tZ$.
Thus $\bigcup_{T_1,T_2} M_{T_1,T_2} = L$, 
so it suffices to show the generalized shape theorem on $M_{T_1,T_2}$. 
By the shape theorem for $Z$ there is a $\mu >0$ 
independent of $\Ga,\tGa,\be$ such that 
\[
(1-\ep/2)B(0,\be\mu^{-1}) \subset Z_t/t \subset (1+\ep/2)B(0,\be\mu^{-1})
\quad 
\text{ for sufficiently large $t$}. 
\]
In order to see that this implies 
\[
(1-\ep)B(0,\be\mu^{-1}) \subset (\tZ_t \cup \Ga)/t \subset
(1+\ep)B(0,\be\mu^{-1}) \quad 
\text{ for sufficiently large $t$} 
\]
on $M_{T_1,T_2}$,  
we use that $\frac{t-T_2}{t} \ge \frac{1-\ep}{1-\ep/2}$ 
for sufficiently large $t$ and 
\begin{equation} \label{shapecompare}
\tZ_{t} \; \cup \; \Ga \subset Z_{t} \;\subset \; 
\tZ_{t+T_2} \cup \Ga  \quad  
\text{ for all } t \ge 0 \quad \text{ on } M_{T_1,T_2}.
\end{equation}
For the proof of \eqref{shapecompare} we fix 
a point configuration $X \in M_{T_1,T_2}$. We note that 
the chosen coupling of $\tZ$ and $Z$ implies 
that every $p = (x,s,r,w) \in X$ with $x \notin \Ga$ 
and $w \le \be$ produces an outburst in both models, 
and in both models the time between the infection of $x$ 
and the time $p$ produces the outburst is $s$. \\
For the inclusion $\tZ_t \cup \Ga \subset Z_t$ we observe that 
$Z_t \supset \Ga$ is trivial and $\tZ_t \subset Z_t$
can be shown by induction on the growth times of $\tZ$: 
For $t=0$ the assertion is trivial. If $t>0$ is a growth 
time, then either this growth is from the initial configuration 
$\tGa$ of $\tZ$ (and in this case $\tZ_t \subset Z_t$ follows
from $\Ga \subset Z_t$) 
or there is a point $p = (x,s,r,w) \in X$ producing 
an outburst at time $t$ w.r.t.\ $\tZ$. In this case $w \le \be$ 
and $x$ was infected in $\tZ$ at time $t' = t-s$. 
By induction hypothesis $\tZ_{t'} \subset Z_{t'}$, 
so $x$ was infected in $Z$ at an earlier time $t'' \le t'$, 
i.e.\ $p$ produces an outburst w.r.t.\ $Z$ at time $t''+s \le t'+s = t$.  
\\
For $Z_{t} \subset \tZ_{t+T_2} \cup \Ga$ we argue 
with induction on the growth times of $Z$: 
For $0 \le t \le T_1$ the assertion is trivial since  
$Z_{T_1} \subset \tZ_{T_2} \cup \Ga$ on $M_{T_1,T_2}$.
If $t>T_1$ is a growth time w.r.t.\ $Z$, there is a  point 
$p=(x,s,r,w)$ producing an effective outburst in $Z$ at time $t$. 
So $x$ was infected in $Z$ at time $t' = t - s< t$, which gives 
$x \in Z_{t'} \subset \tZ_{t'+T_2} \cup \Ga$ by induction hypothesis. 
On $M_{T_1,T_2}$ the last effective outburst located in $\Ga$ occurs before $T_1$. 
Since $t>T_1 $ we know $x \notin \Ga$, so $x\in \tZ_{t'+T_2 } \cup\Ga $ 
implies $x \in \tZ_{t'+T_2}$. 
Thus the outburst caused by $p$ w.r.t.\ $\tZ$ occurs no 
later than $(t'+T_2) + s = t+T_2$.


\section{Proof of the main result: Theorem \ref{thmgro}} \label{secproofthm}

\subsection{Outline of the strategy} \label{secoutline}

We denote the two given \CRMs\ by $Z(\Ga)$ and $Z(\Ga')$. 
Assuming 
\begin{equation} \label{assump}
\pr(\gro(\Ga)) > 0, \;  \quad \text{ and } \quad  
\pr(\esc_i(\Ga')) > 0 \; \text{ for } i=1,2,  
\end{equation}
it suffices to show $\pr(\gro(\Ga')) > 0$.
W.l.o.g.\ we may assume $\be_1 \ge \be_2$. 
Our basic strategy to show $\pr(\gro(\Ga')) > 0$ 
in the next subsections will be to decompose the evolution of $Z(\Ga)$
into an initial part and a final part separated by a space-dependent 
time horizon. This will be chosen such that with positive probability
the initial evolution will infect certain fundamental 
regions before the time horizon, and in the final evolution 
outbursts in these fundamental regions will yield infinite infection paths  
for both types after the time horizon. 
After that we will describe an initial evolution of $Z(\Ga')$ infecting 
the fundamental regions before the time horizon and having positive probability.
By coupling the final evolution of $Z(\Ga')$ to the final evolution 
of $Z(\Ga)$ we also have coexistence for the initial 
condition $\Ga'$, i.e. $\pr(\Ga') >0$.  


\subsection{Escaping the initial configuration in $Z(\Ga')$}
\label{secescaping}

By \eqref{assump} the initial configuration $\Ga'$ allows 
each type to escape some ball with positive probability.
However, what we need is  that $\Ga'$ also allows both types 
to escape some large ball at the same time with positive probability. 
To describe these escape routes we need some notation for line segments: 
For given points $x_1 \neq x_2$ let $\sli_{x_1x_2}$ denote 
the straight line passing through these points. 
We identify $\sli_{x_1x_2}$ with $\R$ such that $x_1 < x_2$ 
so that we can use the induced order on $\sli_{x_1x_2}$ along with 
the corresponding interval notation for line segments. 

\begin{defi} \label{defiexit}
We say that an infection configuration $(A_1,A_2)$ has 
escape corridors specified by $x_1,x_2 \in \R^d$, $r_0,\de_0 >0$
whenever $d(x_1,x_2) > \rp + 4 \de_0$, 
$B(x_i,\de_0) \subset A_i \subset B(0,r_0)$,
$(-\infty,x_1]_{+\de_0} \cap A_2 = \emptyset$
and $[x_2,\infty)_{+\de_0} \cap A_1 = \emptyset$.
\end{defi}
We note that in the situation of Definition \ref{defiexit} we can use 
Lemma~\ref{legrowth} to let the two species grow along $(-\infty,x_1]$ and $[x_2,\infty)$, 
and we will show that this situation occurs with positive probability: 
\begin{lem} \label{leexit}
Let $\rp >0$ be the constant chosen in \eqref{radprorp}. 
For suitably chosen $T_0' \ge 0$, $x_1,x_2 \in \R^d$, 
$\de_0 < {\rp}/ {8}$ and $r_0 >2 \rp$ we have $\pr(M_0')>0$, 
where $M_0'$ is the event that $Z_{T_0'}(\Ga')$ 
has escape corridors specified by $x_1,x_2,r_0,\de_0$. 
\end{lem}


\subsection{Time horizon and fundamental regions for $Z(\Ga)$}
\label{sec:Timehorfundreg}

In the following we will define subsets of $M_0 := \gro(\Ga)$ 
by specifying certain parameters of the infection evolution 
of $Z(\Ga)$. After that we will have to show that the parameters 
can be chosen in a way such that these subsets of $M_0$ 
have positive probability. In Step~1 we use the 
parameters $r_0,T_0'$ chosen in Lemma~\ref{leexit}.
In Step~3 we will restrict our attention to only one infinite infection path
of the stronger type. This idea was already used by M.~Deijfen et al.\ 
in \cite{DH2} for the two-type discrete Richardson-model. 
\\

\noindent{\bf Step 1:}
Let $r_2 > r_1 > r_0$ and $B_i:=B(0,r_i)$. 
Let $M_1$ be the set of configurations of $M_0$ 
such that 
\begin{itemize}
\item 
at time $T_0'$ the infected region is contained in $B_1$,
\item
no outburst with position in $B_1$ infects anything outside of $B_2$. 
\end{itemize}
\noindent{\bf Step 2:}
Let $\die>0$. Let $M_2$ be the set of configurations of $M_1$
such that $2\die$ is a lower bound on the distances of the positions 
of effective outbursts in $B_2$ to 
\begin{itemize}
\item
the lines passing through the position of another effective outburst 
in $B_2$ and the origin,
\item
the boundaries of infection balls of outbursts that infect 
regions in $B_2$, 
\item
the boundaries of the balls $B_1$ and $B_2$.
\end{itemize}
\noindent{\bf Step 3:}
Let  $0<\dio<\die$, $D^1, D^2_1,\ldots,D^2_k$ balls of radius $\dio$ 
(``fundamental balls'') and 
$T(D^1), T(D^2_1),\ldots,T(D^2_k) >0$ times associated to the balls. 
Let $M_3$ be the set of configurations of $M_2$
such that
\begin{itemize}
\item 
$D^1$ contains the starting point 
of an infinite type-1 infection path that never returns to $B_2$,
\item
the effective outbursts of type 2 within $B_2 - B_1$ 
are located in the balls $D^2_1,\ldots,D^2_k$ 
($k \ge 1$, one ball for each outburst),
\item 
at the time $T(.)$ corresponding to a fundamental ball
the position of the corresponding outburst is infected, 
but the outburst has not yet occurred.
\end{itemize}
\begin{lem}\label{lespecify}
The above parameters $r_1$, $r_2$,  $\die$, $\dio$,  
$D^1,\ldots,D^2_k$, $T(D^1), \ldots, T(D^2_k)$ 
can be chosen so that 
$\pr(M_i) > 0$ for $i=1,2,3$, and furthermore
\begin{equation} \label{proppara}
r_1 >  r_0 + 4\rp, \quad 
r_2 > r_1 + 4 \rp, \quad \die < \tfrac \rp 4, \quad \dio < 
\min\{\tfrac \die {10}, \tfrac {\die^2}{8\rp}\}.
\end{equation}
\end{lem}
We now choose the parameters according to the above lemma, and 
consider them to be fixed for the rest of the proof.
$\Mli := M_3$ is the desired specification of the (initial and final)
infection evolution of $Z(\Ga)$. We have $\pr(\Mli) >0$ 
by Lemma~\ref{lespecify}, and on $\Mli$ the following properties hold:
\begin{itemize}
\item 
All fundamental balls are contained in $B_2-B_1$ and 
keep a distance of at least $2\die-\dio$ to the boundaries of this set. 
All type-2 fundamental balls keep a distance of at least
$2\die-\dio$ to the line through the center of $D^1$ and the origin. 
The distance between two fundamental balls is at least 
$2\die-2\dio$. The times corresponding to the fundamental balls 
are $>T_0'$.
\item
Each fundamental ball is infected as a whole (so the 
type-1 fundamental ball is completely 1-infected and 
the type-2 fundamental balls are completely 2-infected),
and each fundamental ball contains the position of exactly 
one effective outburst. A fundamental ball is infected at 
the associated time, but the effective outburst in 
this ball occurs later. 
\item 
Every 2-infection of a region in $B_2^c$ originating in 
$B_2$, in fact originates in one of the type-2 
fundamental balls. The type-1 
fundamental ball is the starting point of an infection path 
to infinity for type 1 such that this path never returns to $B_2$.
\end{itemize}
In Section~\ref{secoutline} we introduced the notion of a space-dependent time horizon.
We now choose the time horizon on the fundamental balls to be 
the times associated with these balls, and 
on $B_2^c$ to be $T(B_2^c):=0$. On the rest of the space it will be chosen later.


\subsection{Initial evolution of $Z(\Ga')$}

In this section we only consider the model with initial 
configuration $\Ga'$, so we set $Z := Z(\Ga')$.
We would like to construct an initial evolution starting from 
$\Ga'$ that infects all fundamental balls before the time horizon. 
More precisely we define 
\[
T' := \min\{T(D^1), T(D^2_1),\ldots,T(D^2_k) \} > T_0'
\]
and $\Mpl$ to be the set of all configurations 
such that for $Z$ at time $T'$ ball $D^1$ is 1-infected, 
all balls $D^2_i$ $(1 \le i \le k)$ are 2-infected, 
and no point of $B_2^c$ is infected.
By definition $\Mpl$ only depends on the point configuration in  
\[ 
\Din :=  B_2 \times [0,T']. 
\]
In order to show $\pr(\Mpl) > 0$ we give a step by step description 
of a suitable infection evolution with the desired properties.
The main reason why this explicit construction is possible 
is that we only have to infect one small ball with type 1, 
but are allowed to infect everything else in $B_2$ with type 2. 
For the construction we choose fixed times 
$T_0' < T_1' < \ldots < T_4' := T'$, and we define a sequence of  
events $M_i'$ describing a desirable infection state at time $T_i'$.\\ 
\begin{figure}[p] 
\begin{center}
\psfrag{A}{$B_0$}
\psfrag{B}{$B_1$}
\psfrag{C}{$B_2$}
\psfrag{D}{$C$}
\psfrag{E}{$C_{+\rp-\dio} -B_1$}
\psfrag{F}{$B(c_2,\dio)$}
\psfrag{G}{$B(c_1,\dio)$}
\psfrag{H}{$C_+$}
\psfrag{I}{$D^1$}
\psfrag{J}{$D^2_j$}
\includegraphics[scale = 0.45]{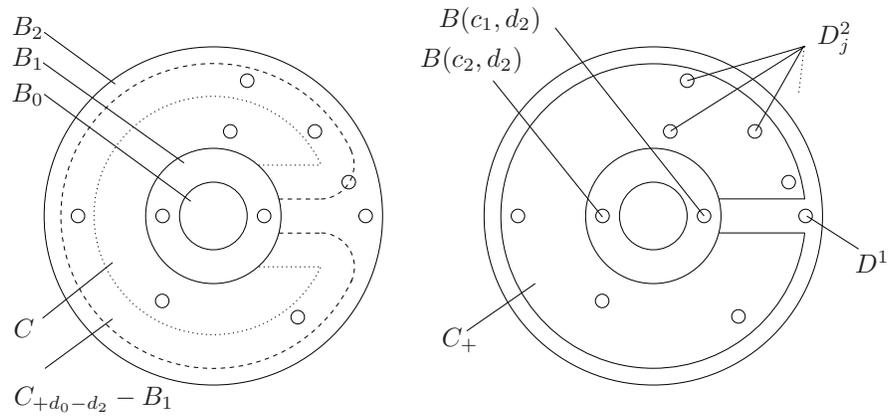}
\end{center} 
\caption{Sketch of the sets relevant in the definition of the initial evolution} 
\label{figinitial}
\end{figure}
\begin{figure}[p] 
\psfrag{A}{$T_1'$}
\psfrag{B}{$T_2'$}
\psfrag{C}{$T_3'$}
\psfrag{D}{$T_4'$}
\psfrag{1}{= infected by species 1}
\psfrag{2}{= infected by species 1 or not at all}
\psfrag{3}{= infected by species 2}
\psfrag{4}{= infected by species 2 or not at all}
\psfrag{5}{= not infected}
\psfrag{6}{= infection state unspecified}
\begin{center}
\includegraphics[scale = 0.45]{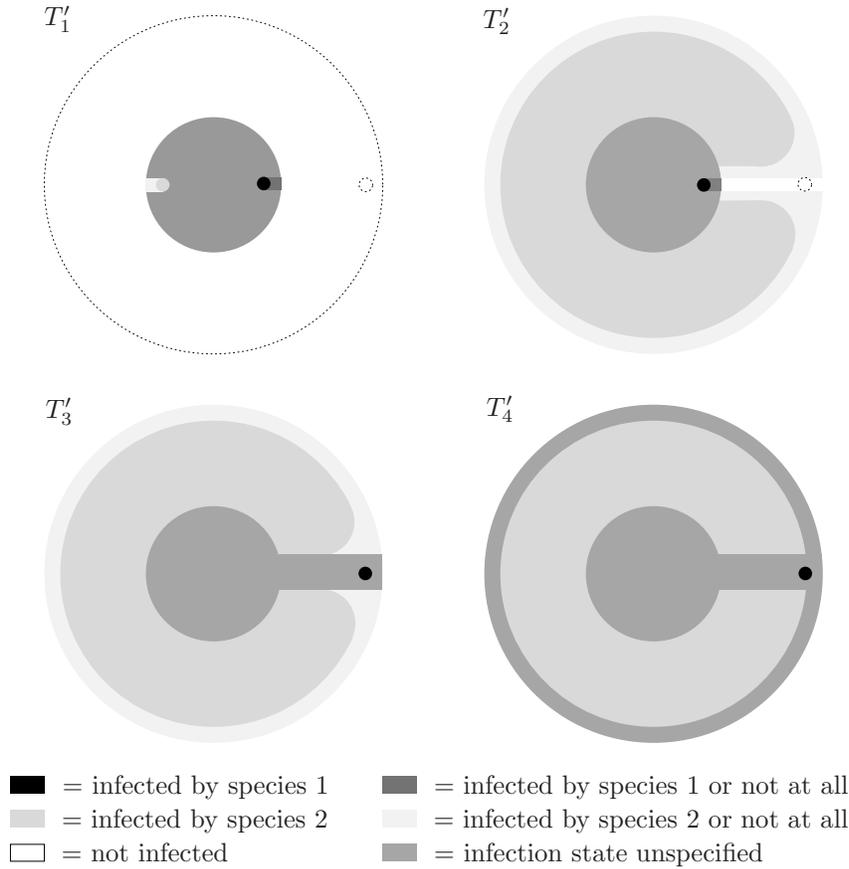}
\end{center} 
\caption{Sketch of the infection states at time $T_i'$ on 
the set $M_i'$ for $1 \le i \le 4$.} 
\label{figinitialevolution}
\end{figure}
We will use interval notation for line segments of $\sli_{0a}$, 
where $a$ denotes the center of $D^1$. 
Furthermore let $c_1$ denote 
the point on $[0,\infty) \subset \sli_{0a}$ with 
$d(c_1,0) = \tfrac{r_0+r_1} 2$,  $c_2 := -c_1$, 
$C := B_2 - B_1 - (\partial B_2)_{+\rp+\dio} - [0,\infty)_{+\rp+2\dio}$ and  
$C_{+} :=  B_2 - B_1 - (\partial B_2)_{+\die} - [0,\infty)_{+\die}$, see Figure \ref{figinitial}.
We note that $D^1\subset [0,\infty )_{+\die }$ and thus $C_{+}\cap D^1=\emptyset $.
The events $M_i'$ (as shown in Figure \ref{figinitialevolution}) can now be 
defined by
\[
\begin{split} 
M_1' &:= \{ B(c_1,\dio) \subset Z_{T_1'}^1 \subset B_1 - [-r_1,c_2]_{+\dio},
B(c_2,\dio) \subset Z_{T_1'}^2 \subset B_1 - [c_1,r_1]_{+\dio}\}\\
M_2' &:= \{ B(c_1,\dio) \subset Z_{T_1'}^1 \subset B_1, 
C_{+\rp-\dio} -B_1 \subset Z_{T_2'}^2 \subset B_2 - 
[c_1,r_2]_{+\dio}\} \\
M_3' &:=\{ D^1 \subset Z_{T_3'}^1 \subset B_1 \cup [0,r_2]_{+\die} \cap B_2,
C_{+\rp-\dio} -B_1 \subset Z_{T_3'}^2 \subset B_2\}\\
M_4' &:= \{ D^1 \subset Z_{T_4'}^1 \subset B_2, 
C_{+} \subset Z_{T_4'}^2 \subset B_2\}.
\end{split}
\]
$M_0'$ is the event of positive probability chosen in Lemma~\ref{leexit}.
\begin{lem} \label{leplppr}
There are events $M'_{(i-1) \to i}$ $(i=1,2,3,4)$ 
describing an infection evolution in $(T_{i-1}',T_i']$ such that 
$\pr(M'_{(i-1) \to i})>0$ and $M_{i-1}' \cap M'_{(i-1) \to i}\subset M_i'$.
\end{lem}
By the independence property of the Poisson process,
$\pr(M_{i-1}' \cap M'_{(i-1) \to i}) = \pr(M_{i-1}')\pr( M'_{(i-1) \to i})$, so Lemma~\ref{leplppr} inductively gives 
$\pr(M_i')>0$ for $i=0,\ldots,4$. As by construction all balls $D^2_i$ 
are  contained in $C_{+}$, we have $M'_4 \subset \Mpl$ and thus 
\begin{equation} \label{mpl}
\pr(\Mpl) > 0.
\end{equation}


\subsection{Final evolution of $Z(\Ga')$} \label{secinf}

Let us now consider the 2-type \CRM\ $\tZ = Z(\tGa)$ with 
generalized initial configuration $\tGa$ consisting of 
all fundamental balls with their associated times and types and 
$D_0$, the remaining part of $B_2$, with associated time $\infty$. 
We assume that $\tZ$ is constructed 
using the Poisson point process underlying $Z(\Ga)$. 
On $\Mli$ the infection evolution of $\tZ$ in $B_2^c$ 
is very similar to the one of $Z(\Ga)$, differing only 
in that type 1 has been reduced considerably. 
However, we will show that in $\tZ$ type 1 is still strong enough 
to grow without bound. Let $\Mi := G(\tGa)$ denote the event that in 
$\tZ$ both types grow without bounds. 
\begin{lem} \label{lereduction}
We have $\Mli \subset \Mi$ and thus $\pr(\Mi) >0$. 
\end{lem}
We now would like to use $\Mi$ for the final evolution of $Z(\Ga')$, 
but unfortunately on $\Mpl \cap \Mi$ we are not guaranteed to have coexistence for $Z(\Ga')$. 
The problem is that in the construction of $\tZ$ from a point configuration 
all points in $D^0$ are ignored, 
so $\Mi$ has no information on points in this region, 
whereas using $\Mpl$ for the initial evolution of $Z(\Ga')$
it is very likely that $D^0$ will be at least partially infected at 
time $T'$. To guarantee coexistence in $Z(\Ga')$, we thus would like to 
delete any point of the underlying point configuration in $D^0$ 
producing a virtual effective outburst virtually infecting 
points of $B_2^c$ after time $T'$. By Lemma~\ref{leeffectivegen} 
we know that $\tleo_{B_2} < \infty$ a.s. on $\Mi$, because 
$\Mi \subset L$. 
Thus we can choose a time $T'' > T'$ such that for 
\begin{equation} \label{mpippr}
\Mpi := \Mi \cap \{\tleo_{B_2} < T''\} \quad 
\text{ we still have } \quad \pr(\Mpi) >0.
\end{equation} 
$\Mpi$ will be used for the final evolution in 
the $\Ga'$-model, and setting $T(D^0) := T''$ we have defined 
a time horizon on all of $\R^d$. Since $\Mpi$ is defined in terms of $\tZ$, 
it only depends on the part of the point configuration after the time horizon, 
i.e. on the point configuration in  
\[ 
\Dfi := \bigcup_{D \in \Dset} D \times [T(D),\infty), 
\quad \text{ where } \quad 
\Dset = \{B_2^c,D^0,D^1,D^2_1,\ldots,D^2_k\}.
\]
We now delete the remaining points that could interfere 
with the desirable final evolution.  Let $\Mpdel$ be the set 
of all point configurations without any point in 
\[
\Dde := \bigcup_{D \in \Dset} D \times [T',T(D)) \times \R_+ \times 
[0,\be_1].
\]
As $T(B_2^c)=0$ we have $\la_{\times}^d(\Dde) < \infty$, 
and thus 
\begin{equation} \label{mpdel}
\pr(\Mpdel) > 0.
\end{equation} 
By construction the sets of point configurations $\Mpl$, $\Mpdel$ 
and $\Mpi$ depend on $\Din$, $\Dde$ and $\Dfi$ respectively, 
and these sets are disjoint. By the independence 
property of the Poisson process we thus have 
\[
\pr(\Mpli) = \pr(\Mpl) \pr(\Mpdel) \pr(\Mpi) >0 \; 
\text{ for } \;
\Mpli := \Mpl \cap \Mpdel \cap \Mpi,  
\]
where all probabilities are positive because of \eqref{mpl}, 
\eqref{mpippr} and \eqref{mpdel}. It is an immediate consequence 
of the construction that on $\Mpli$ the model $Z(\Ga')$ has exactly 
the same infection evolution in $B_2^c$ as the model $\tZ$, and 
we have coexistence for the model $\tZ$ on $\Mpli$
since $\Mpli \subset \Mpi \subset \Mi$.  
This implies that we also have coexistence for the model $Z$ on $\Mpli$,
i.e. $\pr(\gro(\Ga')) >0$.


\section{Proof of the Lemmas} \label{secprooflemthm}

\subsection{Escape corridors: Lemma \ref{leexit}}

In this subsection we only consider the \CRM\ with initial configuration $\Ga'$, 
so let $Z := Z(\Ga')$. $Z$ doesn't change a.s.\ if $\Ga'$ is modified outside 
of its support. Thus we may assume w.l.o.g.\ that 
$\supp(\Ga_i') = \overline{\Ga_i'}$, which implies $\supp(Z_i^t) = \overline{Z_i^t}$
for all $t \ge 0$. 
The proof of Lemma \ref{leexit} consists of three parts. 
First we construct an event $M_0''$ of positive probability 
such that on $M_0''$ there is a time $T$ and there are distinct 
$y_i \in \overline{Z_T^i}$ $(i=1,2)$ such that 
using interval notation on $\sli_{y_1y_2}$
\begin{equation} \label{corridor}
(-\infty,y_1-\rp] \cap \overline{Z_T^2} = \emptyset \;
\text{ and } \;
[y_2+\rp,\infty) \cap \overline{Z_T^1} = \emptyset. 
\end{equation}
In the second part we show that a stronger version of \eqref{corridor}
holds for fixed values of $T, y_1, y_2$ with positive probability.
In the third part we verify the existence of deterministic escape corridors.\\

For the first part we define the minimal enclosing ball (MEB) of a set $K \subset \R^d$ 
to be a closed ball containing $K$ with minimal radius. 
Every bounded set has a unique MEB. 
We will make repeated use of the following property of MEBs: 
\begin{equation} \label{meb}
\text{The MEB $B$ of a compact set $K$ is also the MEB of 
$\partial B \cap K$.}
\end{equation}
Let $B = B(m,r)$ denote the MEB of $\overline{\Ga_1'} \cup \overline{\Ga_2'}$.
W.l.o.g. $\partial B \cap \overline{\Ga_1'} \neq \emptyset$ and  
$\partial B \cap \overline{\Ga_2'} = \emptyset$. 
(By \eqref{meb} one of these sets is nonempty. If it is the other way round
we can argue similarly. If both are nonempty the assertion is trivial 
as we can choose $T=0$ and $y_i \in \partial B \cap \overline{\Ga_i'}$.)
Thus by \eqref{meb} $B$ is the MEB of $\partial B \cap \overline{\Ga_1'}$.\\
{\bf Case 1:} $\rho([r,\infty)) >0$. In this case we may 
assume w.l.o.g.\ that $\rp$  with property \eqref{radprorp} was chosen such that 
$\rp \ge r$.
We set $M_0'' := \X$, $T=0$ and choose $y_2 \in \overline{\Ga'_2} - \{m\}$
and $y_1 \in \partial B \cap \overline{\Ga_1'}$ such that $d(y_1,y_2) >r$. 
(As $B(m,r)$ is the MEB of $\partial B \cap \overline{\Ga_1'}$, 
$B(y_2,r)$ doesn't enclose this set.) 
In order to see that $y_1,y_2$ have the above property, 
we only have to check 
$[y_2+\rp,\infty) \cap \overline{\Ga_1'} \subset [y_2+\rp,\infty) \cap B = \emptyset$, which follows from 
$\la^1(\sli_{y_1y_2} \cap B - [y_1,y_2]) < 2r - r = r \le \rp$.\\
{\bf Case 2:} $\rho([r,\infty)) = 0$. Let $M_1'' := \{Z_{T''}^1 = \Ga_1'\}$ 
and $M_2'' = \{Z_{T''}^2 \not \subset B\}$, where $T''$ is chosen sufficiently large such that $\pr(M_2'')>0$. We also have $\pr(M_1'')>0$, 
and since the absence of outbursts of type 1 can only strengthen type 2
we have $\pr(M_2''|M_1'')>\pr(M_2'')$. 
Setting $M_0'' = M_1'' \cap M_2''$ we obtain $\pr(M_0'') > 0$.
On $M_0''$ at some time $T$ an outburst 
$B' = B(x',r')$ of species 2 occurs that for the first time 
infects a region in $B^c$. As $r' < r$ a.s. and as  
$B$ is the MEB of $\partial B \cap \overline{\Ga_1'}$, 
$B'$ can't enclose this set. So we can choose points
$y_1 \in \partial B \cap \overline{\Ga_1'}-B'$ and $y_2 \in \partial B \cap B'$
that have the above property.\\

For the second part we note that the infection configuration is 
always constant over a time interval 
of positive length, $\overline{Z_T^i}$ is always a compact set, 
the infected region is always bounded, and 
if the points $y_i'$ are sufficiently close to $y_i$, these
points still satisfy \eqref{corridor}. 
This shows that on $M_0''$ there are rational parameters 
$T' >0$, $0 < \de < \rp/8$, $y_1',y_2' \in \R^d$ and $r>0$  
such that  
\begin{equation} \label{precorridor}
\begin{split}
(-\infty,&y_1'-\rp]_{+4\de} \cap \overline{Z_{T'}^2} = \emptyset, \; 
[y_2'+\rp,\infty)_{+4\de}\cap  \overline{Z_{T'}^1} = \emptyset, \; 
d(y_1',y_2')> 4\de,\\
& Z_{T'}^{\cup} \subset B(0,r) \; \text{ and } \;  
\la^d(B(y_i',\de) \cap Z_{T'}^{i})>0 \text{ for } i = 1,2. 
\end{split}
\end{equation}
Thus for some particular choice of these parameters the event $M_{00}'$ 
that \eqref{precorridor} holds for this choice of parameters has positive probability.\\ 

Finally, for the last part let $M_{01}'$ be the set of all point configurations 
with exactly one point in each of the sets
\[
\big( B(y_i',\de) \cap Z_{T'}^{i}\big) \times (T',T'+1] \times [\rp,\rp + \de) \times [0,\be_2]
\quad (i=1,2),
\]
and no additional points in  $B(0,r) \times (T',T'+1] \times [0,\be_1]$.
$M_{01}'\in\F_{\R^d\times (T',T'+1]}$ is independent of 
$M_{00}'\in\F_{\R^d\times [0,T']}$ and has positive probability. 
As a consequence  $M_0' := M_{00}' \cap M_{01}'$ has positive probability.
On $M_0'$ the infection configuration $Z_{T'+1}$ 
has escape corridors specified by the parameters 
$x_1 := y_1' - \rp + 2\de$ 
and $x_2 := y_2' + \rp - 2\de$, $r_0 := r + 2\rp$ and $\de_0 := \de$.


\subsection{Specifying the evolution in the $\Ga$-model: 
Lemma~\ref{lespecify}}

\noindent{\bf Step 1:}
We have 
$M_0 \cap \{Z_{T_0'} \subset B(0,r_1) \} \cap \{\leo_{B_1} < T\} 
\cap \{Z_T \subset B(0,r_2) \} \subset M_1$ for any $T >0$. 
The events $\{Z_{T_0'} \subset B(0,r_1) \}$, $\{\leo_{B_1} < T\}$ and $\{Z_{T} \subset B(0,r_2)\}$ are increasing in $r_1$, $T$ and $r_2$ respectively, and $\leo_{B_1} < \infty$ 
a.s.\ by Lemma~\ref{leeffective} and for fixed $T$ $~Z_T$ is bounded a.s..
Thus we only have to choose first $r_1$, then $T$ and then 
$r_2$ large enough in order to obtain $\pr(M_1) >0$.\\

\noindent{\bf Step 2:}
Observing that the event of all given distances being bounded 
from below by $2\die$ is decreasing in $\die$, we are left to show 
that the infimum of the distances is positive a.s.. 
Since every single distance is positive a.s.,  
it suffices to observe that the number of effective outbursts 
in $B_2$ is finite a.s. by Lemma~\ref{leeffective} and the number 
of outbursts that infect a region in $B_2$ is finite a.s., 
as $B_2$ is completely infected after a finite time 
by Lemma~\ref{leeffective}.  
Therefore if we choose $\die>0$ small enough, we have $\pr(M_2)>0$. \\

\noindent{\bf Step 3:}
We choose a small $\dio>0$ satisfying \eqref{proppara}, 
and we fix a finite covering of $B_2$ consisting of balls 
of radius $\dio$. For any infection configuration it is possible 
to choose a finite number of balls with the properties described
in Subsection \ref{sec:Timehorfundreg}, Step 3, 
 and as in each ball the time span between the infection of the 
outburst position and the outburst at that position is positive 
a.s.\ we can choose a rational time serving as the time horizon 
on the ball. The number of ways to choose a finite set of fundamental 
balls from the finite covering and a rational time horizon is countable, so 
there is a specific choice of fundamental balls and corresponding time horizons 
such that $\pr(M_3)>0$.


\subsection{Local evolution in the $\Ga'$-model: 
Lemma \ref{leplppr}}

\begin{figure}[!htb] 
\begin{center}
\psfrag{a}{$x_1$}
\psfrag{b}{$x_2$}
\psfrag{A}{$C_1$}
\psfrag{B}{$C_2$}
\psfrag{c}{$c_1$}
\psfrag{d}{$c_2$}
\psfrag{h}{$B_0$}
\psfrag{f}{$B_{1/2}$}
\psfrag{x}{$b_1$}
\psfrag{y}{$b_2$}
\includegraphics[scale = 0.45]{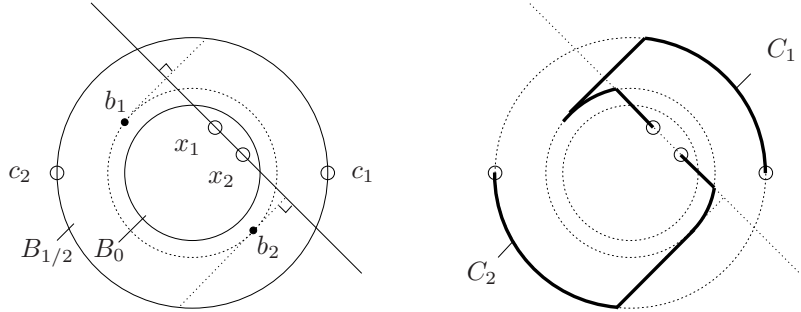}
\end{center} 
\caption{Construction of the curves $C_1$ and $C_2$} 
\label{figcorridor}
\end{figure}

$M_{0\to 1}'$: Infecting the balls $B(c_i,\dio)$.
We use the parameters chosen in Lemma~\ref{leexit},  $\de := \de_0$
and $T_{1/2}' := (T_0'+T_1')/2$.
Choosing two curves 
$C_i \subset B_{1/2} :=  B(0,r_0+2\rp)$ $~(i = 1,2)$ 
joining $x_i \in B_0$ and $c_i \in \partial B_{1/2}$, 
we use Lemma~\ref{legrowth} to let type 1 grow along 
$C_1$ in the time interval $(T_0',T_{1/2}']$ 
with precision $\de$  and after that to let type 2 
grow along $C_2$ in the time interval $(T_{1/2}',T_1']$
with precision $\de$. 
This infection evolution defines an event $M_{0 \to 1}'$ with the desired properties
provided that the curves $C_i$ satisfy 
$d(C_1,Z_{T_0'}^2) > \de$, $d(C_2,Z_{T_0'}^1) > \de$, 
$d(C_1,C_2) = d(x_1,x_2) >  \rp + 4 \de$ and 
$d(C_1,B(c_2,\dio)) > \rp + 2 \de$. 
For a possible construction of such curves see Figure~\ref{figcorridor}. Here $b_1$ and $b_2$ are 
antipodal points on the surface of $B(0,r_0 + \frac{\rp}{2})$ 
such that $\sli_{b_1b_2}$ is parallel to $\sli_{x_1x_2}$.
After passing $b_i$ the curve $C_i$ continues on 
the line in the hyperplane perpendicular to $\sli_{b_1b_2}$ 
that hits $B_{1/2}$ in a point with minimal distance to $c_i$.

\noindent
$M_{1\to 2}'$: Infecting most of $B_2-B_1$ with type 2.
We use Lemma~\ref{legrowth} to let type~2 grow along $C' := C \cup [-r_1,c_2]$ 
in the time interval $(T_1',T_2']$ with precision $\de := \dio/2$. 
This infection evolution defines an event $M_{1 \to 2}'$ with the desired properties.
\newpage 

\begin{figure}[!htb] 
\begin{center}
\psfrag{a}{$c_1$}
\psfrag{b}{$a_0$}
\psfrag{c}{$a_1$}
\psfrag{d}{$a_2$}
\psfrag{e}{$a_3 = a$}
\psfrag{f}{$S_{d_1}(b_0)$}
\psfrag{g}{$S_{d_1}(b_1)$}
\psfrag{h}{$S_{d_1}(b_2)$}
\psfrag{A}{$C_{+\rp-\dio}$}
\psfrag{B}{$\partial B_2$}
\includegraphics[scale = 0.5]{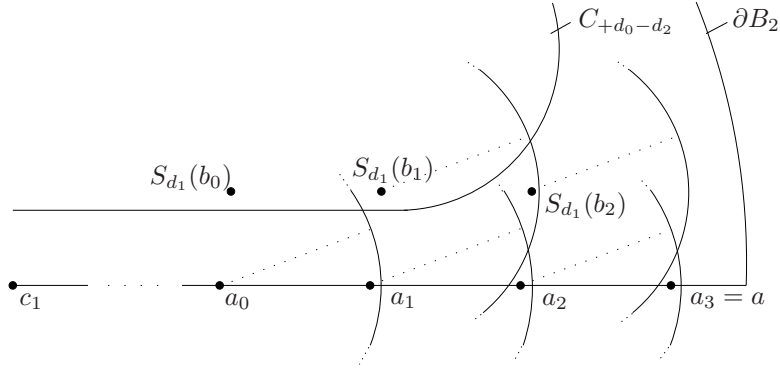}
\end{center} 
\caption{Definition of the infection evolution of step $M_{2\to 3}'$
in the difficult case. 
Some infection outbursts are indicated by segments of their boundaries, which are connected to their centers by dotted lines.
} 
\label{figjump}
\end{figure}

\noindent
$M_{2\to 3}'$: Infecting $D^1$ with type 1.
In this step we basically would like to let type 1 
grow along $[c_1,r_2]$ with precision $\de := \dio/2$ 
until it reaches $D^1$. In order to
keep the 1-infections close to $[c_1,r_2]$, we have to ensure 
that for every outburst of type 1 all points of the corresponding 
ball that are too far away from $[c_1,r_2]$ already are 2-infected 
before the outburst occurs. In this case we say that species 2 provides
sufficient containment for the outburst. In the following
things get more complicated if $D^1$ is near $\partial B_2$, 
and thus outbursts infecting $D^1$ do not yet have 
sufficient containment. If so, we let type 1 grow along $[c_1,r_2]$ 
only as long as there is sufficient containment. After that we 
alternatingly produce outbursts of type 2 (providing further containment
and possibly infecting part of $[c_1,r_2]$) and outbursts of type 1 
(infecting points close to $[c_1,r_2]$ further up), see Figure~\ref{figjump}. 
It turns out that we have to repeat this scheme at most twice. 
In the remainder of this subsection we will give a more detailed description  
of this construction. \\
For a line segment $\sli'$ of the line 
$\sli = \sli_{c_1c_2}$ let $S_r(\sli')$ denote the cylindrical shell 
with axis $\sli'$ and radius $r$, i.e.\ the set of all points $x$ with  
$d(x,\sli) = r$ such that the projection of $x$ onto $\sli$ is 
on $\sli'$. 
In case of a single point $x$ we write $S_{r} (x):=S_{r} (\{x\})$.
In order to determine for which points $a' \in [c_1,a]$ an 
infection outburst at $a'$ with precision $\de$ 
has sufficient containment on $M_2'$, we first note that 
\begin{equation} \label{distrand}
d(S_{r}(r_2 -s),\partial B_2) \ge s-r \quad 
\text{ for all } 0 \le r \le s \le r_2.
\end{equation}
(This and other purely geometric relations 
will be proved at the end of this subsection.) \eqref{distrand} 
implies  
$S_{\rp + 2\dio}([0, r_2 - 2\rp - \die])-B_1  \subset C$
and thus 
\begin{equation} \label{distrandcon}
S_{3\dio}([0, r_2 - 2\rp - \die]) \subset 
S_{\rp + 2\dio}([0, r_2 - 2\rp - \die])_{+\rp-\dio} 
\subset C_{+\rp-\dio}  \cup B_1.
\end{equation}
On $M_2'$ we know that $C_{+\rp-\dio} -B_1$ is 2-infected, 
so this shows that the above described infection outburst at 
$a'$ has sufficient containment if $a' \in [c_1,a]$ is such that 
\begin{equation}\label{tubeinfect}
a' \le r_2 - 2 \rp - \die - (\rp + 2 \de) = 
r_2 - (3\rp +\die+\dio).
\end{equation}
Thus, if $a' := a-(\rp-2\de)$ satisfies \eqref{tubeinfect} 
(i.e. if $a \le r_2 - (2\rp + \die + 2\dio)$) we can use 
Lemma~\ref{legrowth} to let type 1 grow along $[c_1,a']$ 
in the time interval $[T_2',T_3']$ with precision $\de$, 
and immediately get an event $M'_{2 \to 3}$ as desired.
However, if $a >  r_2 - (2\rp + \die +2\dio)$, 
we might need two more outbursts (of size $\rp$) to reach $D^1$ 
with containment not yet provided. 
We now define $a_0 < a_1 < a_2 < a_3 := a$ and $b_i \le a_{i+1}$
to be those points on $[c_1,a]$ satisfying  
\begin{equation} \label{grdefab}
 d(a_i,a_{i+1}) = \rp-2\dio \quad \text{ and } \quad  
d(S_{\die}(b_i),a_{i+1}) = \rp + 2\dio. 
\end{equation}
We note 
\begin{equation} \label{grlageb}
 a_i \le b_i \le a_i + \die/2 \quad \text{ and } \quad a_3 = a \le r_2 - 2\die . 
\end{equation}
Figure~\ref{figjump} shows a halfplane starting from $\sli_{c_1a}$. 
The three points above $[c_1,a]$ indicate the points of intersection of the halfplane 
with the cylindrical shells $S_{\die}(b_i)$. We would like to define 
$M_{2\to 3}'$ to be the following infection evolution (where  
$T_2' = t_0' < t_1' < \ldots < t_5' = T_3'$): 
\begin{enumerate}
\item growth of type 1 along $[c_1,a_0]$ in $[t_0',t_1']$, 
\item growth of type 2 along $S_{\die}([b_0,b_1])$ in $[t_1',t_2']$,
\item growth of type 1 by a single outburst at $a_1$ 
in $[t_2',t_3']$,
\item growth of type 2 along $S_{\die}([b_1,b_2])$ in $[t_3',t_4']$,
\item growth of type 1 by a single outburst at $a_2$ 
in $[t_4',t_5']$.
\end{enumerate}
In every step the growth is meant to be within $B_2$ and with precision $\de$.
The important part of the boundary of the infection outburst in 
every step is indicated in Figure~\ref{figjump}.
By Lemma~\ref{legrowth} we have $\pr(M_{2\to 3}')>0$, and 
Figure~\ref{figjump} shows that for $M_2' \cap M_{2\to 3}' \subset M_3'$ 
we only have to ensure the following properties in the 
corresponding steps: In Step~1 we need sufficient containment, 
i.e.\ $a_0$ has to satisfy \eqref{tubeinfect}. 
This is a consequence of \eqref{grdefab}, \eqref{grlageb} and $8 \dio \le \die$. 
Also, $B(a_1,\dio)$ has to be fully infected, 
which follows from \eqref{grdefab}. 
For Step~2 we observe that $b_0 \le r_2 - 2\rp - \die$, 
which follows from \eqref{grdefab} and \eqref{grlageb}. 
Thus \eqref{distrandcon} implies $S_{\die}(b_0)_{+\dio} \subset C_{+\rp-\dio}$, 
which means that the growth of type 2 can start from $S_{\die}(b_0)$. 
No part of $B(a_2,\dio)$ is 2-infected by \eqref{grdefab}.
Also Step~2 has to provide sufficient containment for 
the next, i.e.\ we need 
\begin{equation} \label{grnoinfa}
 B(a_i,\rp+2\dio) \subset [0,r_2]_{+\die} \cup 
S_{\die}([b_{i-1} ,b_i])_{+(\rp-\dio)} 
\end{equation}
for $i=1$. In Step~3 $B(a_2,\dio)$ has to be fully infected, 
which follows from the first part of \eqref{grdefab} 
Steps 4 and 5 have the analogous properties of Steps 2 and 3, 
provided we have \eqref{grnoinfa} for $i=2$, too. Finally, to ensure 
that the growth is always within $B_2$ we need 
\begin{equation}\label{grcontain}
S_{\die}([b_0,b_2])_{+(\rp+\dio)} \cup  B(a_2,\rp+2\dio)
\subset B_2. 
\end{equation}
What remains to be checked are properties \eqref{distrand}, \eqref{grlageb}, \eqref{grnoinfa} and  \eqref{grcontain}.
\begin{figure}[htb!]
\begin{center}
\psfrag{A}{\eqref{distrand}}
\psfrag{B}{\eqref{grlageb}}
\psfrag{a}{$r_2-s$}
\psfrag{b}{$s$}
\psfrag{c}{$r$}
\psfrag{d}{$S_r(r_2-s)$}
\psfrag{e}{$\partial B_2$}
\psfrag{1}{$a_i$}
\psfrag{3}{$a_{i+1}$}
\psfrag{4}{$S_{d_1}(b_i)$}
\psfrag{5}{$\rp-2\dio$}
\psfrag{6}{$\die$}
\psfrag{7}{$\rp+2\dio$}
\includegraphics[scale = 0.35]{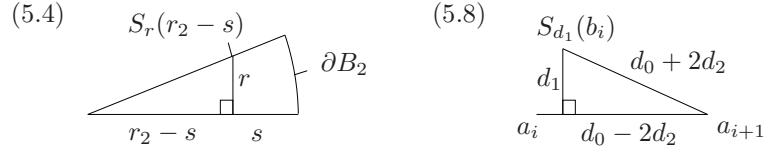}
\caption{Illustrations for proofs of \eqref{distrand} and \eqref{grlageb}}
\label{figdistrand}
\end{center} 
\end{figure}
As can be seen from Figure~\ref{figdistrand},
assertion \eqref{distrand} is equivalent to 
\[
\sqrt{(r_2-s)^2 + r^2} + (s-r) \le r_2,
\]
which holds for all $r \le s \le r_2$. Figure~\ref{figdistrand}
also shows that the first part of \eqref{grlageb} is equivalent to 
\[
0 \; \le \; \rp - 2 \dio - \sqrt{(\rp+2\dio)^2 - \die^2} \; 
\le \; \die/2, 
\]
which is satisfied for $\dio \le \die^2/(8\rp)$ and $2 \die \le \rp$. 
The second part of \eqref{grlageb} follows by definition of $a$ and $\die$. 
\begin{figure}[htb!]
\begin{center}
\psfrag{A}{\eqref{distrand}}
\psfrag{B}{\eqref{grlageb}}
\psfrag{a}{$a_i$}
\psfrag{b}{$S_{d_1}(b_i)$}
\psfrag{c}{$\rp+2\dio$}
\psfrag{d}{$\rp-\dio$}
\psfrag{e}{$S_{d_1}(b_i +\rp - \dio)$}
\psfrag{f}{$a_{i+1}$}
\psfrag{3}{$a_{i+1}$}
\psfrag{4}{$S_{d_1}(b_i)$}
\psfrag{5}{$\rp-2\dio$}
\psfrag{6}{$\die$}
\psfrag{7}{$\rp+2\dio$}
\includegraphics[scale = 0.45]{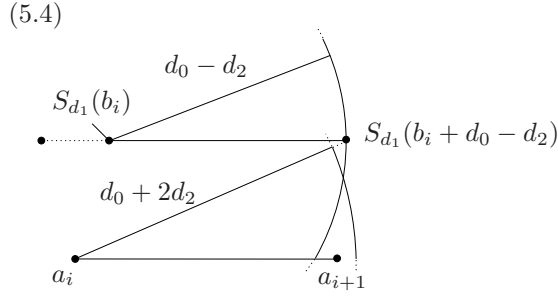}
\caption{Illustration for proof of \eqref{grnoinfa}}
\label{figcontainment}
\end{center} 
\end{figure}
For \eqref{grnoinfa} we observe that 
$[a_i - (b_i-a_i),a_i + (b_i-a_i)] \subset [b_{i-1},b_{i}]$,
which follows from \eqref{grlageb},
so by symmetry it suffices to check that 
\[
d(S_{\die}(b_i + \rp -\dio),a_i) 
\ge d(S_{\die}(b_i), a_{i+1}) = \rp + 2 \dio,
\]
see Figure~\ref{figcontainment}. Here the equality follows from 
\eqref{grdefab} and Figure~\ref{figcontainment} shows that the 
inequality follows from  $a_i \le b_i$ and $b_i + \rp - \dio \ge 
a_i + \rp - \dio > a_{i+1}$. For \eqref{grcontain} we first 
observe that $ B(a_2,\rp+2\dio) \subset B_2$
follows from $r_2 - a_2 \ge 2 \die + (\rp - 2\dio) \ge \rp+2\dio$, 
where we have used \eqref{grlageb}. Furthermore we have 
\[
S_{\die+\dio}(r_2 - \rp - \die - 2\dio)_{+\rp+\dio} \subset B_2, 
\]
by \eqref{distrand}, so it suffices to show that 
$b_2 \le r_2 - \rp - \die - 2\dio$. This follows from 
\eqref{grdefab} and \eqref{grlageb}.\\

\noindent
$M'_{3 \to 4}$: Infecting the remaining part of $B_2-B_1$ 
with type 2. 
Let $r_2' := d(S_{\die + \dio}(r_2 - \rp - \die - 2\dio),0)$, 
and $C' := B(0,r_2') -B_1 - [0,r_2]_{+\die+\dio}$. 
We define $M'_{3 \to 4}$ by the growth of type 2 along $C'$ 
in the time interval $[T_3',T_4']$ with precision $\de := \dio/2$. By
Lemma~\ref{legrowth} we have $\pr(M'_{3 \to 4}) >0$ and 
for $M'_3 \cap M'_{3 \to 4} \subset M_4'$ it suffices to observe 
that on $M'_3$ no point of $C'_{+\dio}$ is 1-infected, most of $C'$
is 2-infected, by \eqref{distrand} we have $C'_{+\rp+\dio} \subset B_2$, 
and by choice of $r_2'$ we have $C'_{+\rp-\dio} \supset C_{+}$.


\subsection{Restricting the strong type: Lemma \ref{lereduction}}

In this section we restrict our attention to the infection evolutions 
of $Z := Z(\Ga)$ and $\tZ$ for a fixed point configuration 
$X \in \Mli$. By definition of $\Mli$, we know that $Z^2$ grows
without bound and $X$ contains an infinite type-1 infection path $(p_n)_{n\ge 1}$ w.r.t $Z$ such that for $p_1 = (x_1,t_1,r_1,b_1)$ 
we have $x_1 \in D^1$ and $t_1 \ge T(D^1)$. It suffices to show 
that 
\begin{itemize}
\item $\tZ^2_t \cap B_2^c \supset Z^2_t \cap B_2^c$ 
for all $t \ge 0$, and 
\item $(p_n)_n$ is an infinite type-1 
infection path w.r.t $\tZ$.
\end{itemize}
To prove the first assertion we show that 
\[
\forall t \ge 0: \quad 
\tZc_t \subset \Zc_t, \tZ_t^1 \subset Z_t^1, \text{ and }  
Z^2_t \cap B_2^c \subset \tZ^2_t \cap B_2^c
\]
by induction on the times of growth of both models. 
At $t=0$ the above is trivial, so suppose that $t >0$ is a time  
of growth. If the growth is due to the initial configuration of $\tZ$
the inclusions remain true, because  
in the model $Z$ every fundamental ball is infected at its time horizon.
Otherwise the growth is due to an outburst $p = (x,t,r,w) \in X$ in one of the models. 
If $p$ produces an effective type-1 outburst in $\tZ$, the same outburst happens in $Z$. 
If $p$ produces an effective type-2 outburst in $\tZ$, 
there is also an outburst of type-1 or type-2 in $Z$. 
So if $p$ produces an effective outburst w.r.t.\ $\tZ$, all inclusions remain true. 
Now we take a look at the remaining case, 
i.e.\ $p$ produces an effective outburst in $Z$ and not in $\tZ$. 
We only have to worry about the third inclusion, and thus only have to consider the case 
that $p$ produces an effective type-2 outburst in $Z$ such that $B(x,r) \not \subset B_2$
and $x \in B_2$. (If $B(x,r) \subset B_2$ the outburst contributes to neither side, 
and if $x \notin B_2$ the outburst contributes to both sides.) 
In this case $x$ is in one of the fundamental balls and $t$ is 
after the corresponding time horizon. 
Thus the initial configuration of $\tZ$ ensures that 
the same outburst happens in $\tZ$ and the third inclusion remains true. \\

For the second assertion we show by induction that $p_n$ 1-infects $x_{n+1}$
in $\tZ$. By induction hypothesis $x_n$ gets 1-infected in both models at time 
$t_n$, so we have a type-1 outburst at time $t_{n+1}$ in both models. 
Since $\tZc_{t} \subset \Zc_{t}$ at all times, $x_{n+1}$ can't get infected 
in $\tZ$ earlier than in $Z$. Thus the outburst 1-infects $x_{n+1}$ in $\tZ$ aswell. \\

{\bf Acknowledgments:} We would like to thank H.-O.~Georgii 
for bringing this topic to our attention and for 
helpful discussions and remarks. 


\renewcommand{\thesection}{}

\setlength{\parindent}{0cm}

\end{document}